\documentclass[UKenglish, final]{article}
\usepackage[utf8]{inputenc}

\usepackage{amsmath,mathscinet}
\usepackage{amssymb}
\usepackage{todonotes}
\usepackage{graphicx}
\usepackage{mathrsfs}
\usepackage{lipsum}
\usepackage{ragged2e}
\usepackage{xcolor}
\usepackage{proof}
\usepackage{enumitem}

\usepackage{babel, csquotes, graphicx, textcomp, kantlipsum, booktabs}
\usepackage{hyperref}
\hypersetup{colorlinks,allcolors=black}
\usepackage{varioref, cleveref} 

\RequirePackage{amsfonts}       
\RequirePackage{amsmath}        
\RequirePackage{amssymb}        
\RequirePackage{amsthm}         
\RequirePackage{cancel}         
\RequirePackage{comment}        
\RequirePackage{dsfont}         
\RequirePackage{enumitem}       
\RequirePackage{etoolbox}       
\RequirePackage{mathrsfs}       
\RequirePackage{mathtools}      
\RequirePackage{multirow}       
\RequirePackage{pgffor}         
\RequirePackage{physics}        
\RequirePackage{showkeys}       
\RequirePackage{stmaryrd}       
\RequirePackage{tablefootnote}  
\RequirePackage{textcomp}       
\RequirePackage{thmtools}       
\RequirePackage{tikz}           

\declaretheorem[style = plain, numberwithin = section]{theorem}
\declaretheorem[style = plain,      sibling = theorem]{corollary}
\declaretheorem[style = plain,      sibling = theorem]{lemma}
\declaretheorem[style = plain,      sibling = theorem]{proposition}
\declaretheorem[style = definition,      sibling = theorem]{observation}

\declaretheorem[style = definition, sibling = theorem]{example}

\declaretheorem[style = remark,     sibling = theorem]{remark}

\crefname{observation}{observation}{observations}

\DeclareMathOperator\Aut{Aut}

\newcommand{\g}{\phi}

\title{The Fermat curves, arrangements of lines, and intersections of osculating curves}
\author{Torgunn Karoline Moe and Nils Peder Astrup Toft}

\begin{document}

\maketitle

\begin{abstract}
    In this paper we present new results about arrangements of lines and osculating curves associated to the Fermat curves in the projective plane. We first consider the sextactic points on the Fermat curves and show that they are distributed on three grids. The grid lines constitute new line arrangements and examples of free curves associated with the Fermat curves. Moreover, we compute the hyperosculating conics to the Fermat curves, study the arrangement of these conics, and find that they intersect in a special way. The latter result is a consequence of the action of the group of automorphisms on osculating curves, and we conclude with a more general result for intersections of osculating curves of any given degree.
\end{abstract}

\tableofcontents

\section{Introduction}
Arrangements of lines, conics and curves of low degree in the projective plane has been a popular research topic the last decade. This work has been driven by the interest in the so-called \emph{free} (and \emph{nearly free}) curves, where a key ingredient is the study of singularities and intersections of the irreducible parts of an arrangement. We refer to the work by Dimca in \cite{Dim17}; see also \cite{Dim24} for background and open problems on this topic. 

It turns out that it is possible to construct many examples of free (and nearly free) curves using curves of low degree together with their inflection tangents or hyperosculating conics, lines and conics that have unexpectedly high intersection multiplicity with the curve at a point. Such arrangements have been studied for, amongst others, the nodal cubic, the Klein quartic, and the Fermat curves of degree $d=3$ and $d=4$; see the work by Szemberg and Szpond in \cite{SS24}, Abe, Dimca and Sticlaru in \cite{ADP24}, Dimca, Ilardi, Malara and Pokora in \cite{DGMP25}, and Merta and Zieli{n}ski in \cite{MZ24}. In addition, there are some recent results for the Fermat curve of degree $d \geq 3$ and arrangements of inflection tangents given by Dimca, Ilardi, Pokora and Sticlaru in \cite{DIPS24}. 

With this backdrop, the purpose of the present article is to study the Fermat curves $\mathcal{F}_d$ of degree $d \geq 3$ and arrangements of some associated curves. We do this to get an overview of the situation, in particular with respect to the osculating and hyperosculating curves. A returning theme is how symmetries control the geometry. 

In \Cref{sec:fermat} we recall well-known results about the Fermat curve $\mathcal{F}_d$ and its group of automorphisms. Moreover, we take a closer look at the osculating curves to a point on a Fermat curve, and we use Cayley's formula to compute the osculating conic to any given point on a Fermat curve. In particular, we study the points for which there exists a hyperosculating line or conic, that is, inflection points and sextactic points, and we compute the hyperosculating conics. 

In \Cref{sec:lines} we find a couple of interesting line arrangements associated with the Fermat curves. We show that the sextactic points are clustered in three grids and that arrangements of lines from these grids provide new examples of free curves.

In \Cref{sec:hyposc} we consider intersections of certain osculating curves and show that for the Fermat curves these intersections are a bit surprising. Most noteworthy, our study of the sextactic points and their osculating curves can be summarized in the following result.
\begin{theorem}\label{thm:main}
    Given any set of $d$ co-linear sextactic points on a Fermat curve $\mathcal{F}_d$ of degree $d>3$, the intersection of the corresponding $d$ tangent lines contains exactly one point; a $d$-fold point in the union of the tangent lines. Moreover, there are exactly two points contained in the intersection of the $d$ hyperosculating conics, so the union contains two $d$-fold points.
\end{theorem}

The above result is a consequence of the action of the group of automorphisms on osculating curves, and we end this paper by generalizing our results to osculating curves of higher degree.

\section{The Fermat curves}\label{sec:fermat}
Let $\mathbb{P}^2$ denote the projective plane over the complex numbers $\mathbb{C}$, and let $S=\mathbb{C}[x,y,z]$ be the graded polynomial ring. The plane Fermat curve $\mathcal{F}=\mathcal{F}_d$ of degree $d \geq 3$ is given as the zero set $V(F)$ of the homogeneous polynomial $F=F_d \in\mathbb{C}[x,y,z]_d$, where \[F_d=x^d+y^d+z^d.\] The Fermat curves are all smooth curves since the partial derivatives never vanish at the same point, and the genus is \[g=\frac{(d-1)(d-2)}{2}.\]

\subsection{The group of automorphisms}
We denote the group of automorphisms of the Fermat curve by $G_d$. For $d \geq 4$, it is shown in \cite{Tze95} that the group has order $6d^2$ and is generated by $S_3$ and $\mu(d)\times\mu(d)$, where $\mu(d)$ is the group of roots of unity in $\mathbb{C}$ generated by the primitive root $\zeta$. We choose $\rho,\varphi,\psi:\mathbb{P}^2\to\mathbb{P}^2$ as generators of $G_d$, where
\begin{align*}
   \rho(x:y:z)&=(\zeta x:y:z),\\
    \varphi(x:y:z)&=(y:x:z),\\
    \psi(x:y:z)&=(z:y:x).     
\end{align*}
For $d=3$, the group of automorphisms is well known with similar properties; see f.ex. \cite{DGMP25}.

A crucial observation for our results in \Cref{sec:hyposc} is that some of the $6d^2$ automorphisms fix a line pointwise; i.e., the automorphism restricts to the identity on the line. In fact, the chosen generators have this property; $\rho$ fixes the line $x=0$, $\varphi$ the line $x=y$, and $\psi$ the line $x=z$.

\subsection{Osculating and hyperosculating curves}
Given a plane curve $\mathcal{C}$ and a smooth point $p \in \mathcal{C}$, there exists a unique curve $\Sigma_p$ of degree $n$ such that the intersection multiplicity $(\mathcal{C}. \Sigma_p)_p \geq \frac{1}{2}n(n+3)$. The curve $\Sigma_p$ is referred to as \emph{the osculating curve of degree $n$} to $\mathcal{C}$ at $p$. When the inequality is strict, $\Sigma_p$ is called \emph{hyperosculating}. 

In this section we study the osculating and hyperosculating lines and conics of the Fermat curves.

\subsubsection{The tangent lines}
Given a Fermat curve $\mathcal{F}$ and a point $p$, the tangent line $T_p$ is the unique line for which $(T_p.\mathcal{F})_p\geq2$. Calculating this line is classical and straightforward.
\begin{lemma}\label{lem:tan}
    The tangent line $T_p$ to the Fermat curve $\mathcal{F}_d$ of degree $d$ at a point $p=(p_x:p_y:p_z)$ is given by
\[T_p \colon p_x^{d-1}x+p_y^{d-1}y+p_z^{d-1}z=0.\]
\end{lemma}

\subsubsection{The inflection points and inflection tangents}
It is well known that the Fermat curve $\mathcal{F}_d$ has $3d$ inflection points, points $p$ where the tangent line $T_p$ intersects $\mathcal{F}$ with intersection multiplicity $(T_p.\mathcal{F})_p\geq 3$, see for example \cite{Wat99}. 

This fact can be verified by intersecting $\mathcal{F}$ with the Hessian curve $H(\mathcal{F})=V(H)$, where
\begin{align*}
    H=d^3(d-1)^3(xyz)^{d-2}.
\end{align*}
It follows that the coordinates of the $3d$ inflection points are 
\begin{equation*}\label{thm:inflcoord}
    \{(0:1:u^k),(u^k:0:1),(1:u^k:0): \text{$k$ odd}\},
\end{equation*}
where $u\in\mathbb{C}$ is a $2d$-root of unity, so that $u^d=-1$. 

Moreover, all inflection points are maximal in the sense that the intersection multiplicity $(T_p.\mathcal{F})_p=d$, which can be computed directly using one of the inflection points, say $p=(0:1:u^k)$, with tangent line $T_p$ given by
\begin{align*}
    T_p \colon &y-u^{-k}z=0,
\end{align*}
which has no further intersections with $\mathcal{F}$. This is sufficient since any two inflection points can be related by an automorphism of $\mathcal{F}$. 

\subsubsection{The osculating conics}
The osculating conic $O_p$ to the Fermat curve $\mathcal{F}$ at a point $p$, for which $(O_p.\mathcal{F})_p \geq 5$, can be computed using Cayley's formula from 1859 \cite{Cay59}. Indeed, after careful computations, the defining polynomial for the osculating conic at a point on a Fermat curve turns out to be quite simple.
\begin{lemma}\label{lem:osccon}
    The osculating conic $O_p$ to a point $p=(p_x:p_y:p_z)$ on the Fermat curve $\mathcal{F}_d$ of degree $d$ is given by the polynomial 
   \small{
    \begin{align*}
        &p_x^{2d-2}(d+1)\left((2d-1)p_y^dp_z^d+(2-d)(p_y^d+p_z^d)p_x^d\right)x^2\\
        &+p_y^{2d-2}(d+1)\left((2d-1)p_x^dp_z^d+(2-d)(p_x^d+p_z^d)p_y^d\right)y^2\\
        &+p_z^{2d-2}(d+1)\left((2d-1)p_x^dp_y^d+(2-d)(p_x^d+p_y^d)p_z^d\right)z^2\\
        &-\left(2(d+1)(d-2)p_x^{2d-1}p_y^{2d-1}+4(2d-1)(d-2)(p_x^{d-1}p_y^{2d-1}+p_x^{2d-1}p_y^{d-1})p_z^d\right)xy\\
        &-\left(2(d+1)(d-2)p_x^{2d-1}p_z^{2d-1}+4(2d-1)(d-2)(p_x^{d-1}p_z^{2d-1}+p_x^{2d-1}p_z^{d-1})p_y^d\right)xz\\
        &-\left(2(d+1)(d-2)p_y^{2d-1}p_z^{2d-1}+4(2d-1)(d-2)(p_y^{d-1}p_z^{2d-1}+p_y^{2d-1}p_z^{d-1})p_x^d\right)yz.
    \end{align*}}
\end{lemma}

\begin{proof}
   With notation as in \cite{MM19}, by Cayley's formula from \cite{Cay59}, $O_p$ is given by the polynomial 
    \begin{equation}\label{eq:osccon}
        9H^3(p)D^2F_p-\left(6H^2(p)DH_p+9H^3(p)\Lambda(p)DF_p\right)DF_p,
    \end{equation}
    where $9H^3\Lambda=-3\Omega H+4\Psi$. 
    
    We obtain the following expressions for $\Omega$ and $\Psi$ via straightforward calculations:
    \begin{align*}
        \Omega&=d^5(d-1)^5(d-2)(d-3)(xyz)^{d-4}(x^dy^d+y^dz^d+z^dx^d),\\
        \Psi&=d^8(d-1)^8(d-2)^2(xyz)^{2d-6}(x^dy^d+y^dz^d+z^dx^d).
    \end{align*}
    Substituting back into our formula for $9H^3\Lambda$ leads to
    \begin{align*}
        9H^3(p)\Lambda(p)&=-3\Omega(p) H(p)+4\Psi(p)\\
        &=d^8(d-1)^8(d+1)(d-2)(p_xp_yp_z)^{2d-6}(p_x^dp_y^d+p_y^dp_z^d+p_z^dp_x^d).
    \end{align*}
  Moreover, we find the following evaluations and polynomials
    \begin{align*}
        H^2(p)&= d^6(d-1)^6(p_xp_yp_z)^{2d-4},\\
        H^3(p)&= d^9(d-1)^9(p_xp_yp_z)^{3d-6},\\
        DF_p&= d(p_x^{d-1}x+p_y^{d-1}y+p_z^{d-1}z),\\
        D^2F_p&= d(d-1)(p_x^{d-2}x^2+p_y^{d-2}y^2+p_z^{d-2}z^2),\\
        DH_p&= d^3(d-1)^3(d-2)(p_xp_yp_z)^{d-3}(p_yp_zx+p_xp_zy+p_xp_yz).
    \end{align*}
    The result follows by substituting these expressions into Cayley's formula for the osculating conic (\ref{eq:osccon}). 
\end{proof}

\subsubsection{The sextactic points and hyperosculating conics}
To identify the points on $\mathcal{F}$ where the osculating conic is hyperosculating, the so-called sextactic points, we intersect the curve with the classical $2$-Hessian curve $H_2(\mathcal{F})$. 

The defining polynomial $H_2=H_2(F)$ of the $2$-Hessian for the Fermat curve can be computed using the formula by Cayley from \cite{Cay65}, corrected in \cite{MM19}, which in this situation reduces to
\begin{align*}
    H_2&=
    \begin{vmatrix}
x^{5d-9} & x^{4d-9} & x^{3d-9}\\
y^{5d-9} & y^{4d-9} & y^{3d-9}\\
z^{5d-9} & z^{4d-9} & z^{3d-9}\\ 
\end{vmatrix}\\
    &=(xyz)^{3d-9}(x^d-y^d)(y^d-z^d)(z^d-x^d).
\end{align*}

Removing the factors of $H_2$ that are factors in $H$, and thus intersects $\mathcal{F}$ at the (higher order) inflection points, we are left with a core that intersects $\mathcal{F}$ at the sextactic points only.

This core can be factorized further, with $\zeta$ a primitive $d$-root of unity, to \[\prod_{j=0}^{d-1}\left(x-\zeta^jy\right)\prod_{j=0}^{d-1}\left(y-\zeta^jz\right)\prod_{j=0}^{d-1}\left(z-\zeta^jx\right).\]
We denote this arrangement of lines by $\mathcal{B}$, and write \[\mathcal{B} \colon \left(x^d-y^d\right)\left(y^d-z^d\right)\left(z^d-x^d\right)=0,\] and we let \[\mathcal{B}_z=V\left(x^d-y^d\right), \quad \mathcal{B}_x=V\left(y^d-z^d\right), \quad \mathcal{B}_y=V\left(z^d-x^d\right).\]

By intersecting the lines $\mathcal{B}$ with the Fermat curve $\mathcal{F}$, we compute the coordinates of its sextactic points. The lines $B_z=V(x-\zeta^jy)$ intersected with $\mathcal{F}$ for all $j \in \{0,\ldots, d-1\}$ and $k$ odd with $0< k< 2d$, gives the $d^2$ points \[s_{j,k}=(\zeta^{j}:1:u^{-k}2^{1/d}).\] The sextactic points on $\mathcal{B}_y$ and $\mathcal{B}_x$ can be found by permutation of the coordinates, and there are in total exactly $3d^2$ sextactic points on $\mathcal{F}$. Observe that given any two sextactic points there is a symmetry of $\mathcal{F}$ relating them. 

Since $\mathcal{F}$ is smooth, we can apply the formula in \cite[Theorem~1.2]{MM19} to show that the $3d^2$ hyperosculating conics $O_s$ at the sextactic points all have the property that $(O_s,\mathcal{F})_s=6$. Indeed, since $(T_p,\mathcal{F})_p=d$ for all inflection points $p$, we get from the formula that the number of sextactic points (counted with multiplicity) is equal to \[6(2d+5g-5)-3d(4+4d-15)=3d^2.\] Since we computed $3d^2$ sextactic points above, the formula ensures that there are no sextactic points for which $(O_s,\mathcal{F})_s>6$.

\begin{remark}
    These points are exactly the so-called \emph{Leopoldt Weierstrass points}, see \cite{Roh82}.
\end{remark}

The defining polynomials of the hyperosculating conics to the Fermat curve can be computed directly from \Cref{lem:osccon}.
We do this for the $d^2$ sextactic points on $\mathcal{B}_z$. 
\begin{corollary}\label{cor:hyposccon}
    Let $s_{j,k}:=(\zeta^j:1:u^{-k}2^{1/d})$ be a sextactic point on $\mathcal{B}_z$. Then the hyperosculating conic to $s_{j,k}$ is given by the polynomial
    \begin{align*}
        O_{j,k}=&d(d+1)\zeta^{-2j}x^2\\
        &+d(d+1)y^2\\
        &-4(d+1)(2d-3)u^{2k}2^{-2/d}z^2\\
        &-2(d-2)(5d-3)\zeta^{-j}xy\\
        &+8d(d-2)\zeta^{-j}u^{k}2^{-1/d}xz\\
        &+8d(d-2)u^{k}2^{-1/d}yz.
    \end{align*}
\end{corollary}

\section{Line arrangements associated with sextactic points}\label{sec:lines}
With $\mathcal{C} \colon P=0$ a reduced curve of degree $\hat{d}$ in $\mathbb{P}^2$ and $J_P=(P_x,P_y,P_z)$ the Jacobian ideal of $P$, let \[AR(P)=\left\{(a,b,c) \in S^3 \colon aP_x+bP_y+cP_z=0 \right\}\] denote the graded $S$-module of Jacobian syzygies of $P$. By \cite{Dim17}, a reduced plane curve $\mathcal{C}$ is said to be \emph{free} if $AR(P)$ is minimally generated by two homogeneous syzygies of degrees $r$ and $\hat{d}-1-r$, referred to as the exponents. Moreover, by \cite[Corollary~1.2]{Dim17} the curve $\mathcal{C}$ is free if and only if $r \leq (\hat{d}-1)/2$ and
    \begin{equation}\label{eq:freeness}r^2-(\hat{d}-1)r+(\hat{d}-1)^2=\tau(V(P)),\end{equation} where $\tau({V(P)})$ is the total Tjurina number of the arrangement $V(P)$.

There are several line arrangements that are associated with the Fermat curves. First of all, the union of all the tangent lines at inflection points gives the arrangement $\mathcal{A}$, where 
\[\mathcal{A}\colon (x^d+y^d)(y^d+z^d)(z^d+x^d)=0,\] which in different combinations with $\mathcal{F}$ and $\mathcal{H}$ provides examples of free curves, see \cite[Theorem~1.6 and Theorem~1.7]{Dim24}.

In this section we study arrangements of curves starting out with the $3d$ lines $\mathcal{B}$. In particular, we find that the sextactic points are distributed on three grids, where each grid contains $3d$ lines. The lines in these grids can be arranged to provide new examples of free curves.

\subsection{The \texorpdfstring{$\mathcal{B}$}{B}-lines}
We first give a description of the $\mathcal{B}$-lines.

\begin{proposition}\label{thm:2hint}
     The line arrangement $\mathcal{B}$ contains $d^2$ triple points and three ordinary $d$-fold points at the origins. Neither of these points are contained in $\mathcal{F}$.
\end{proposition}

\begin{proof}
    The $d$ lines in $\mathcal{B}_z=V(x^d-y^d)$ all contain $(0:0:1)$, so it is a $d$-fold point. By symmetry, the same hold for $(0:1:0)$ and $(1:0:0)$. 

    Moreover, taking any line $x-\zeta^{j_1} y=0$ from $\mathcal{B}_z$ and any line $y-\zeta^{j_2}z=0$ from $\mathcal{B}_x$, they intersect in the point $(\zeta^{j_1+j_2}:\zeta^{j_2}:1)$. This point is also contained in the line $z-\zeta^{-j}x=0$ for $j=j_1+j_2$ from $\mathcal{B}_y$, which makes it a triple point in $\mathcal{B}$. There are $d^2$ ways to choose two lines from different $\mathcal{B}_i$, hence $d^2$ such points. 

It can be directly verified that neither the origins, nor any other intersection points, are contained in $\mathcal{F}$.
\end{proof}

\begin{remark}\label{rem:Bfree}
The line arrangement $V(xyz) \cup \mathcal{B}$ is by \cite[Corollary~2.9]{MV23} known to be a free curve with exponents $(d+1,2d+1)$. Similar results hold for subarrangements by \cite[Corollary~2.10]{MV23}; in particular, $\mathcal{B}$ is free with exponents $(d+1,2d-2)$.    
\end{remark}

\begin{observation}\label{obs:Bfree}
 The arrangement $\mathcal{F}\cup \mathcal{B}$ of degree $4d$ is \emph{not} free. By \Cref{thm:2hint}, the arrangement has $d^2$ triple points and three $d$-fold points, as well as $3d^2$ double points at the sextactic points. At an ordinary $d$-fold point the Tjurina number is $(d-1)^2$, so the freeness criteria in \Cref{eq:freeness} gives 
 \[r^2-(4d-1)r+(4d-1)^2=d^2(3-1)^2+3(d-1)^2+3d^2.\] This simplifies to \[r^2-(4d-1)r+6d^2-2d-2=0,\] a quadratic polynomial in $r$ with negative discriminant for positive values of $d$, suggesting that $r$ is complex, which is impossible.
\end{observation}

\subsection{The three grids of sextactic points}
The $2$-Hessian by Cayley is not the only curve that intersects a given curve in its sextactic points. Indeed, adding any term (of appropriate degree) with $F$ as a factor to the defining polynomial $H_2$ gives a new curve with similar properties. In particular, and more surprisingly, in the case of the Fermat curve we can work with the factors of degree $d$ of the core of the $2$-Hessian, which ultimately leads to the conclusion that the sextactic points are arranged on three grids.

Before we state the result, let $\mathcal{M}$ and $\mathcal{N}$ denote the following two line arrangements,
\begin{align*}
\mathcal{M} \colon&   \left(z^d+2y^d\right)\left(x^d+2z^d\right)\left(y^d+2x^d\right)=0,\\
& \prod_{\underset{k\text{ odd}}{k=1}}^{2d-1}\left(z-u^{-k}2^{1/d}y\right)\prod_{\underset{k\text{ odd}}{k=1}}^{2d-1}\left(x-u^{-k}2^{1/d}z\right)
    \prod_{\underset{k\text{ odd}}{k=1}}^{2d-1}\left(y-u^{-k}2^{1/d}x\right)=0,
\end{align*}
and
\begin{align*}
\mathcal{N} \colon & \left(y^d+2z^d\right)\left(z^d+2x^d\right)\left(x^d+2y^d\right)=0,\\
& \prod_{\underset{k\text{ odd}}{k=1}}^{2d-1}\left(y-u^{-k}2^{1/d}z\right)\prod_{\underset{k\text{ odd}}{k=1}}^{2d-1}\left(z-u^{-k}2^{1/d}x\right)
    \prod_{\underset{k\text{ odd}}{k=1}}^{2d-1}\left(x-u^{-k}2^{1/d}y\right)=0.
\end{align*}

\noindent As in the case of $\mathcal{B}$, we use the convention that $\mathcal{M}_x=V(z^d+2y^d)$ etc., where the subscript indicates the variable that is \emph{not} in the defining polynomial.

\begin{proposition}
    \label{pro:coremod}
    The $3d^2$ sextactic points on the Fermat curve $\mathcal{F}$ of degree $d\geq 3$ are organized in three clusters of $d^2$ points. 
    
    With $j$ and $k$ as above, the $d^2$ points \[s_{j,k}=(\zeta^j:1:u^{-k}2^{1/d})\]
    are distributed with $d$ points on each of the $3d$ lines in $\mathcal{B}_z=V\left(x^d-y^d\right)$, $\mathcal{M}_x=V\left(z^d+2y^d\right)$, and $\mathcal{N}_y=V\left(z^d+2x^d\right)$. 
    
    The intersection of the $3d$ lines consists of $d^2$ ordinary triple points at the sextactic points and three ordinary $d$-fold points at the origins.   
\end{proposition}

\begin{proof}
Firstly, $\mathcal{B}_z$, $\mathcal{M}_x$, and $\mathcal{N}_y$ are given by bivariate homogeneous polynomials of degree $d$, hence each defining polynomial factors into $d$ linear polynomials.
Secondly, 
    \begin{align*}
    z^d+2y^d&=F-\left(x^d-y^d\right),\\
    z^d+2x^d&=F+\left(x^d-y^d\right).
\end{align*}
Thus, by construction, both $V\left(z^d+2y^d\right)$ and $V\left(z^d+2x^d\right)$ intersect $\mathcal{F}$ in the same $d^2$ sextactic points as $V\left(x^d-y^d\right)$. This implies that through any sextactic point there is exactly one line from each group, and on each of the lines there must be $d$ sextactic points. 

Lastly, the latter claim follows since the $d$ lines in $\mathcal{B}_z$ all contain $(0:0:1)$, the lines in $\mathcal{M}_x$ contain $(1:0:0)$, and the lines in $\mathcal{N}_y$ contain $(0:1:0)$.
\end{proof}

\begin{remark}\label{rem:BMN}
Any union of $\mathcal{B}_i$, $\mathcal{M}_j$ and $\mathcal{N}_k$ where $i \neq j \neq k$, intersected with the Fermat curve gives exactly one cluster of $d^2$ sextactic points. 

Moreover, both the union $\mathcal{B}_i \cup \mathcal{M}_i \cup \mathcal{N}_i$ for $i$ any of $x,y,z$, and any arrangement of lines $\mathcal{K}=\mathcal{K}_x \cup \mathcal{K}_y \cup \mathcal{K}_z$, for $\mathcal{K}$ any of $\mathcal{B}, \mathcal{M},\mathcal{N}$, provide curves that intersect the Fermat curve in all of its $3d^2$ sextactic points. 
\end{remark}

\begin{observation}\label{obs:lines}
   When $d >3$, the only lines containing $d$ sextactic points are the lines in $\mathcal{B}$, $\mathcal{M}$ and $\mathcal{N}$. Indeed, clustering the points in a natural way and writing $\zeta=u^2$, for $j$ and $k$ as above,
   \begin{align*}
   V_z&=s^z(j,k)=\{(u^{2j}:1:u^{-k}2^{1/d})\},\\
   V_y&=s^y(j,k)=\{(1:u^{-k}2^{1/d}:u^{2j}\}),\\
   V_x&=s_x(j,k)=\{(u^{-k}2^{1/d}:u^{2j}:1)\},
   \end{align*}
there is a line containing at least three of them if and only if the determinant of the $(3\times3)$-matrix with the points as rows vanishes.

Up to symmetry, there are three different cases.
\begin{description}[style=nextline, labelwidth=0pt, leftmargin=0pt, labelindent=0pt, align=left]
\item[\;Case I)\;Three points from $V_z$]
    The determinant vanishes if and only if 
    \begin{equation}\label{eq:det3z}
     u^{2j_1-k_3}+u^{2j_2-k_1}+u^{2j_3-k_2}=u^{2j_1-k_2}+u^{2j_2-k_3}+u^{2j_3-k_1}.       
    \end{equation}
 Assume first that $j_1=j_2$. Then either $j_3=j_1=j_2$, and the points are on a line in $\mathcal{B}_z$, or $k_1=k_2$, but then two points are equal. Assume next that $k_1=k_2$. Then either $k_3=k_1=k_2$, in which case the points are on a line in $\mathcal{M}_x$, or $j_1=j_2$, which again is impossible. Lastly, assume that all $j_i$ and $k_i$ are different. Then \Cref{eq:det3z} holds only if $k_i+2j_i\equiv l \pmod{2d} $ for $l$ odd, in which case the points are on a line in $\mathcal{N}_y$.
\item[\;Case II)\;Two points from $V_z$ and one from $V_y$:]
    The determinant vanishes if and only if 
      \begin{equation*}\label{eq:det2z}
     2^{2/d}(u^{2j_2-k_1-k_y}-u^{2j_1-k_2-k_y})+2^{1/d}(u^{-k_2}-u^{-k_1})+(u^{2(j_1-j_y)}-u^{2(j_2-j_y)})=0,       
    \end{equation*}
    which holds only if $j_1=j_2$ and $k_1=k_2$, but then we have chosen the same point twice.
\item[\;Case III)\;One point from each cluster:]
    The determinant vanishes if and only if 
    \begin{equation}\label{eq:detxyz}
          2^{1/d}\left(u^{2j_z-k_y}+u^{2j_y-k_x}+u^{2j_x-k_z} \right)=1+u^{2(j_x+j_y+j_z)}+2^{3/d}u^{-(k_x+k_y+k_z)}+.
    \end{equation}
    Notice that \Cref{eq:detxyz} has a solution if and only if $d=3$; with 
   \begin{align}
         j_x+j_y+j_z&\equiv 0 \pmod{3} \notag\\
          k_x+k_y+k_z &\equiv 3 \pmod{6} \notag
             \end{align}
             
             and
             \begin{equation}\label{detchoose} 
\begin{aligned}
          2(j_y-j_z)-(k_x-k_y)&\equiv 2 \pmod{6}\\
          2(j_x-j_z)-(k_z-k_y)&\equiv 4 \pmod{6}, 
    \end{aligned}  
    \end{equation}
  or any permutation of the variables in (\ref{detchoose}). This can be simplified further by writing $k_i=2m_i+1$. Then there is a solution if, for $j_i,m_i \in\{0,1,2\}$, \begin{align}
         j_x+j_y+j_z&\equiv 0 \pmod{3} \notag\\
          m_x+m_y+m_z &\equiv 0 \pmod{3} \notag
             \end{align}
and any two of the relations
                \begin{equation}\label{detchoosejm} 
\begin{aligned}
          j_x-m_z&\equiv 2 \pmod{3},\\
          j_y-m_x&\equiv 1 \pmod{3},\\
          j_z-m_y&\equiv 0 \pmod{3}.\\
    \end{aligned}  
    \end{equation}
Each ordered triple $\{\left((j_x,m_x),(j_y,m_y),(j_z,m_z)\right) \in (\mathbb{Z}_3 \times \mathbb{Z}_3)^3\}$ corresponds to a line containing three sextactic points. In total, there are nine ordered triples that satisfy the relations; each with six permissible permutations, so altogether 54 lines.
\end{description}
\end{observation}    

\begin{remark}
        When $d=3$, it is shown in \cite[Lemma~4.8]{SS24} that through any sextactic point there are nine lines, each containing three sextactic points. These lines are exactly the ones described in \Cref{obs:lines}; in total 81, with 27 from lines with three sextactic points in the same cluster, and 54 lines with three sextactic points in different clusters. 
\end{remark}

Some of the line arrangements discussed above, and permutations of these (see \Cref{rem:BMN}), turn out to be new examples of free curves.

\begin{theorem}
      The line arrangement $\mathcal{B}_z \cup \mathcal{M}_x \cup \mathcal{N}_y$ is a free curve of degree $3d$ with exponents $(d+1,2d-2)$. Moreover, $V(xyz)\cup\mathcal{B}_z \cup \mathcal{M}_x \cup \mathcal{N}_y$ is a free curve of degree $3d+3$ with exponents $(d+1,2d+1)$. Lastly, the arrangement $\mathcal{F} \cup \mathcal{B}_z \cup \mathcal{M}_x \cup \mathcal{N}_y$ is a free curve of degree $4d$ with exponents $(2d-2,2d+1)$.
\end{theorem}

\begin{proof}
    The first claim can be shown by direct computation of the Jacobian syzygy of the defining polynomial of the line arrangement
\[\mathcal{B}_z\mathcal{M}_x\mathcal{N}_y \colon (x^d-y^d)(z^d+2y^d)(z^d+2x^d)=0,\]   
 which is generated by \[
     \left\{2x^{d+1}-4xy^d+2xz^d, \; -4x^dy+2y^{d+1}+2yz^d, \; -4^{d+1}x^dz-4^{d+1}y^dz-z^{d+1}\right\},\]
     \[\left\{2x^{d-1}y^{d-1},\;-y^{d-1}z^{d-1},\;-x^{d-1}z^{d-1}\right\}.
\]

    By considering the syzygy-module, we see that the exponents are $(d+1,2d-2)$, so $r=d+1$. Moreover, the arrangement is clearly of degree $\hat{d}=3d$.
    
    The left hand side of \Cref{eq:freeness} is then 
    \begin{equation*}
        (d+1)^2-(3d-1)(d+1)+(3d-1)^2=7d^2-6d+3.
    \end{equation*}

    From \Cref{pro:coremod}, we have that the line arrangement has $d^2$ triple points and three $d$-fold points. At an ordinary $d$-fold point the Tjurina number is $(d-1)^2$, so at a triple point, the Tjurina number is $(3-1)^2$. Therefore, the total Tjurina number $\tau(\mathcal{B}_z\mathcal{M}_x\mathcal{N}_y)$ can be found as 
    \begin{equation*}
    d^2 (3-1)^2+3(d-1)^2=7d^2-6d+3
    \end{equation*}
and \Cref{eq:freeness} holds.

The second claim is shown similarly. The arrangement has degree $\hat{d}=3d+3$, we have $r=d+1$, and the left hand side of \Cref{eq:freeness} is then
    \begin{equation*}
        (d+1)^2-(3d+2)(d+1)+(3d+2)^2=7d^2 + 9d + 3.
    \end{equation*}

In the arrangement, there are $d^2$ triple points at the sextactic points, $3d$ double points where each line from $\mathcal{B}_z \cup \mathcal{M}_x \cup \mathcal{N}_y$ intersect $V(xyz)$, as well as three $d+2$-fold points at the origins. Then the right hand side of \Cref{eq:freeness} reduces to \[3(d+1)^2+d^2(3-1)^2+3d=7d^2 + 9d + 3.\]

The third claim can again be shown by direct computation of the Jacobian syzygy of the defining polynomial
\[(x^d+y^d+z^d)(x^d-y^d)(z^d+2y^d)(z^d+2x^d),\]   
 which is generated by 
 \begin{align*}
   \big\{&2x^{2d + 1} - 6xy^{2d} + 6x^{d+1}z^d - 6xy^dz^d + 3xz^{2d}, \\
   &-6x^{2d}y + 2y^{2d + 1} - 6x^dyz^d + 6y^{d+1}z^d + 3yz^{2d},\\
   &-6x^{2d}z - 6y^{2d}z - 6x^dz^{d+1} - 6y^dz^{d+1} - z^{2d + 1}\big\}, \\
   \big\{&-y^{d-1}z^{d-1},\; -x^{d-1}z^{d-1},\; 2x^{d-1}y^{d-1}\big\}.
 \end{align*}

The arrangement has degree $\hat{d}=4d$, we have $r=2d+1$, and the left hand side of \Cref{eq:freeness} is then
    \begin{equation*}
        (2d+1)^2-(4d-1)(2d+1)+(4d-1)^2=12d^2 - 6d + 3.
    \end{equation*}

In the arrangement, there are $d^2$ quadruple points at the sextactic points, in addition to the three $d$-fold points at the origins, so the total Tjurina number sums to
    \begin{equation*}
    d^2 (4-1)^2+3(d-1)^2=12d^2 - 6d + 3,
    \end{equation*}
and \Cref{eq:freeness} holds.

\end{proof}

\subsection{The \texorpdfstring{$\mathcal{M}$}{M}- and \texorpdfstring{$\mathcal{N}$}{N}-lines}
We now study the grid lines and their intersections. We only consider $\mathcal{M}$, as $\mathcal{N}$ (by symmetry) has exactly the same properties.

\begin{proposition}
     The $3d$ lines in $\mathcal{M}$ intersect in $3d^2$ ordinary double points and three ordinary $d$-fold points at the origins. Neither of these points are contained in $\mathcal{F}$.
\end{proposition}

 \begin{proof}
    The $d$ lines in $\mathcal{M}_x=V(z^d+2y^d)$ all intersect in $(1:0:0)$, and similarly for the other two groups and origins, making them $d$-fold points. 

    Moreover, we claim that taking any two lines from different groups, they intersect in a point $q$, which is neither contained in $\mathcal{F}$, nor any line from the third group. By symmetry it suffices to demonstrate this for two lines from $\mathcal{M}_x$ and $\mathcal{M}_{z}$, respectively, say $z-u^{-k_1}2^{1/d}y=0$ and $y-u^{-k_2}2^{1/d}x=0$. These two lines intersect in the point {$\left(1:u^{-k_2}2^{1/d}:u^{-(k_1+k_2)}2^{2/d}\right)$}. By inspection, it proves impossible that this point is contained in any line from $\mathcal{M}_y$, i.e., on the form $x-u^{-k}2^{1/d}z=0$, and similarly it is not contained in $\mathcal{F}$. Thus, the point is an ordinary double point in $\mathcal{M}$.
    
    There are $3d^2$ ways to choose two different lines from three groups, each containing $d$ lines, and the result follows.
\end{proof}

\begin{observation}
     Note that the arrangement $\mathcal{M}$ is \emph{not} free, cf. \Cref{rem:Bfree}. Indeed, the arrangement of degree $3d$ has three ordinary $d$-fold points as well as $3d^2$ ordinary double points. The freeness criteria from \Cref{eq:freeness} then demands that
    \[r^2-(3d-1)r+(3d-1)^2=3(d-1)^2+3d^2(2-1)^2.\] This can be simplified to the quadratic polynomial 
    \[r^2-(3d-1)r+3d^2-2=0,\]
    which has negative discriminant for positive values of $d$, suggesting that $r$ is complex, and this gives a contradiction.  

    Moreover, a similar result can be shown for the union $\mathcal{M} \cup V(xyz)$, again cf. \Cref{rem:Bfree}. Indeed, the arrangement of degree $3d+3$ has three ordinary $(d+2)$-fold points as well as $3d^2$ ordinary double points from intersections in $\mathcal{M}$, and an additional $3d$ double points from intersections between $\mathcal{M}$ and the fundamental triangle. Then \Cref{eq:freeness} demands that
    \[r^2-(3d+2)r+(3d+2)^2=3(d+2-1)^2+3d^2(2-1)^2+3d(2-1)^2.\] This can be simplified to the quadratic polynomial 
    \[r^2-(3d+2)r+3d^2+3d+1=0,\]
    which again has negative discriminant for positive values of $d$.

   Lastly, consider the union $\mathcal{M} \cup \mathcal{F}$, which again is \emph{not} free, cf. \Cref{obs:Bfree}. The arrangement of degree $4d$ has three ordinary $d$-fold points at the origins, as well as $3d^2$ double points in $\mathcal{M}$ and $3d^2$ double points at the sextactic points. Then \Cref{eq:freeness} demands that
    \[r^2-(4d+1)r+(4d-1)^2=3(d-1)^2+3d^2(2-1)^2+3d^2(2-1)^2.\] This can be simplified to the quadratic polynomial 
    \[r^2-(4d-1)r+7d^2-2d-2=0,\]
    which again has negative discriminant.
    
    \end{observation}

\section{Intersections of osculating curves}\label{sec:hyposc}
In this section we take a closer look at the intersections of osculating curves to the Fermat curves. 

\subsection{Fixed lines and intersections of osculating curves}
Due to lack of a suitable reference, we start by showing the invariance of intersections of osculating curves to a plane curve that has automorphisms that fixes lines pointwise.
\begin{theorem} \label{thm:invariantintersection}
    Let $\mathcal{C}\subset \mathbb{P}^2$ be a projective curve invariant under an automorphism $\g\in \Aut(\mathbb{P}^2)$, and let $L\subset\mathbb{P}^2$ be a line fixed pointwise by $\g$. Suppose $P, Q\in \mathcal{C}$ are smooth points in the same orbit under $\g$, and let $\Gamma_P$ and $\Gamma_Q$ be the osculating curves of degree $n$ at $P$ and $Q$, respectively. Then 
    \begin{align*}
        \Gamma_P\cap L=\Gamma_Q\cap L. 
    \end{align*} 
\end{theorem}
\begin{proof}
    Since $P$ and $Q$ are in the same orbit under $\g$, we may choose $k\in \mathbb{Z}$ so that $Q=\g^k(P)$. Note also that since $L$ is a line fixed pointwise under $\g$, it is also fixed pointwise under $\g^k$. Furthermore, since $\g$ is an automorphism and preserves the osculating curves, we have \begin{align*}
        \g^k(\Gamma_P)=\Gamma_{\g^k(P)}=\Gamma_Q.
    \end{align*}
    Putting this together, we get
    \begin{align*}
        \Gamma_P\cap L=\g^k(\Gamma_P\cap L)=\g^k(\Gamma_P)\cap\g^k(L)=\Gamma_Q\cap L,
    \end{align*}
    which completes the proof. 
\end{proof}

Since many automorphisms of the Fermat curve fix lines pointwise, there are a number of geometrical consequences for osculating curves to the Fermat curve that follow from this result. The lines fixed pointwise by an automorphism of $\mathcal{F}$ are exactly the lines of the 2-Hessian. We will consider a few examples in the following sections, with the general result stated in \Vref{cor:fullgenerality}. 

We start with a simple example investigating intersections of tangent lines.

\begin{example}
Let $p$ be any point on $\mathcal{F}$, let $\g$ be an automorphism of $\mathcal{F}$ that fixes a line $L_{\g}$ pointwise, and let $P=\{p_i\}$ be the points in the orbit of $p$ under $\g$. Then the tangent lines $T_i$ to $\mathcal{F}$ at $p_i$ intersect in one point $q \in L_{\g}$.   

In particular, two points $p=(\hat{x}:\hat{y}:1)$ and $q=(\hat{y}:\hat{x}:1)$ on the Fermat curve are in the same orbit under the automorphism $\varphi$. Since $\varphi$ fixes the line $x=y$ pointwise, the tangent lines $T_p$ and $T_q$ intersect in $(1:1:\hat{z})$, where $\hat{z}=-(\hat{x}^{d-1}+\hat{y}^{d-1})$.

Moreover, the $d$ points $p_j=(\zeta^j\hat{x}:\hat{y}:1)$ for $j \in \{0, \ldots, d-1\}$, that are actually co-linear on the line $\hat{y}z-y=0$, are in the same orbit under $\rho$. Since $\rho$ fixes the line $x=0$ pointwise, the $d$ tangent lines to $\mathcal{F}$ at $p_j$ intersect in one point on $x=0$; the point $(0:1:-\hat{y}^{d-1})$.
\end{example}

\subsection{Intersections of tangent lines at inflection points and sextactic points}
Inflection points and sextactic points on a Fermat curve are particularly intriguing, and we inspect the intersections of the tangent lines at these points in depth.

As a second example of an application of \Cref{thm:invariantintersection}, we consider the tangent lines at inflection points.

\begin{example}
  Take the two inflection points $p_x=(0:u^k:1)$ and $p_y=(u^k:0:1)$ for a given $k$. Note that $\varphi(p_x)=p_y$, where $\varphi$ interchanges $x$ and $y$. This automorphism fixes the line $x=y$ pointwise, so $T_{p_x} \colon -u^{-k}y+z=0$ intersects $x=y$ in the same point as $T_{p_y} \colon -u^{-k}x+z=0$, i.e., the point $(1:1:u^{-k})$, which is also easily seen by direct computation.  
\end{example}

A third example shows that $d$ inflection tangents at $d$ co-linear inflection points intersect in only one point. 

\begin{example}
    Consider the $d$ inflection points $p_k=(u^k:1:0)$ on $z=0$ for $k$ odd $0 < k<2d$. These points have tangent lines $T_{p_k} \colon -u^{-k}x+y=0$. The $d$ points are in the same orbit under $\rho$, as $\rho(p_k)=p_{k+2}$. Since $\rho$ fixes the line $x=0$ pointwise, the intersection of $T_{p_k}$ and $x=0$, that is $(0:0:1)$, is fixed. Hence all tangent lines at the $d$ inflection points intersect at this point. This is of course easy to see by direct calculation as well.
\end{example}

Similar results can also be found for the tangent lines at sextactic points. 

\begin{example}
    \label{ex:tangsexthes}
    Consider $d$ sextactic points co-linear on the same line in $\mathcal{B}$, say $V(y-z)$. The points $s_k=(u^{-k}2^{1/d}:1:1)$ are in the same orbit under $\rho$, that fixes the line $x=0$. Moreover, a tangent line \[T_{s_k} \colon \left(u^{-k}2^{1/d}\right)^{d-1}x+y+z=0\] intersects $x=0$ in $(0:-1:1)$, and by \Cref{thm:invariantintersection} (or easy calculation) all the $d$ tangent lines contain $(0:-1:1)$. 
  
Observe that when $d$ is even, $(0:-1:1)$ is not on $\mathcal{F}$, while when $d$ is odd, $(0:-1:1)$ is an inflection point on $\mathcal{F}$.
\end{example}

\begin{remark}
    The latter observation is well known for all smooth cubic curves, including the Fermat cubic, see \cite[Theorem~B]{MZ25}, as the sextactic points correspond to $6$-torsion points that are not $3$-torsion points.
\end{remark}

We have a nearly identical result for groups of $d$ sextactic points co-linear with respect to the $\mathcal{M}$- and $\mathcal{N}$-lines.

\begin{example}\label{ex:tangsextmod}
   Consider $d$ sextactic points co-linear on the same line in $\mathcal{M}$, say one of the lines in $V(z^d+2y^d)$.

 At the $d$ sextactic points $s_j$ with coordinates $s_j=(\zeta^j:1:u^{-k}2^{1/d})$ for $j\in\{0,\ldots,d-1\}$ and $k$ fixed, the curve $\mathcal{F}$ has tangent lines $T_{s_j} \colon \left(\zeta^j\right)^{d-1}x+y+\left(u^{-k}2^{1/d}\right)^{d-1}z=0$. The points are in the same orbit under $\rho$, and $\rho$ fixes the line $x=0$, so the intersection $(0:-u^k2^{(d-1)/d}:1)$ is preserved, thus contained in all the $d$ tangent lines.
  
Lastly, note that $(0:-u^k2^{(d-1)/d}:1)$ is not on $\mathcal{F}$.
\end{example}

\begin{remark}
    Observe that similar results applies to bi-tangents of $\mathcal{F}$. Bi-tangents $T_{q,\hat{q}}$ are lines that are tangents to $\mathcal{F}$ at two points, and it is well known that $\mathcal{F}$ has exactly $\tfrac12d^2(d-2)(d-3)$ such lines \cite{Oka05}. The property of being bi-tangent is preserved by the automorphisms, so by \Cref{thm:invariantintersection} tangents to points in the same orbit under an automorphism will intersect in one point, and this point will be on the line fixed pointwise by the automorphism - if such a line exists.
\end{remark}

\subsection{Intersections of hyperosculating conics}
In this section we go into the intersections of hyperosculating conics both via \Cref{thm:invariantintersection} and through direct computation.

Since a sextactic point is sent to another sextactic point and hyperosculating conics are preserved under the automorphisms, we have the following result as a first corollary of \Cref{thm:invariantintersection}.
\begin{corollary}\label{cor:hyposcgeneral}
 Assume that $g$ is an automorphism of $\mathcal{F}$ with a line $L_{g}$ fixed pointwise by $g$. Let $P_{g}$ be the set of sextactic points in an orbit of $g$. Then the hyperosculating conics at the points in $P_g$ intersect $L_g$ in the same (up to) two points.
\end{corollary}

Taking a more direct approach, we now prove this result for intersections of hyperosculating conics at certain orbits of sextactic points using the defining polynomials of the hyperosculating conics from \Cref{cor:hyposccon}. The advantage of this approach is that we also get the coordinates of the intersection points. For $d=3$, our computation gives a second proof of the result in \cite{DGMP25}.

\begin{theorem}\label{thm:hyposches}
    The $d$ hyperosculating conics to the Fermat curve $\mathcal{F}$ at $d$ sextactic points co-linear on the lines $\mathcal{B}$ have two common intersection points when $d > 3$, and one common intersection point when $d=3$. 
\end{theorem}

\begin{proof}
   Given an integer $j$, let $s_{j,k}:=(\zeta^j:1:u^{-k}2^{1/d})$ for $k$ odd and $0 < k < 2d$ be the $d$ sextactic points on the line $V(x-\zeta^jy)$ in $\mathcal{B}_z$. Let $O_{j,k}$ be the hyperosculating conic to the Fermat curve at the point $s_{j,k}$. 
   
Observe that the $d$ points $s_{j,k}$ are in the same orbit under $g=\psi \circ \rho \circ \psi$, which sends $(x:y:z)$ to $(x:y:\zeta z)$ and fixes the line $z=0$ pointwise. Thus, the hyperosculating conics must intersect $V(z)$ in the same points.

    The intersection of $V(z)$ with $O_{j,k}$ from \Cref{cor:hyposccon} then gives 
    \begin{align*}
        O_{j,k}\cap V(z)&=V\left(d(d+1)\zeta^{-j}x^2-2(d-2)(5d-3)xy+d(d+1)\zeta^jy^2,z\right).
    \end{align*}
    Observe that when $d>3$, the discriminant of the quadratic polynomial,
        \begin{equation*}
    \left(2(d-2)(5d-3)\right)^2-4d(d+1)\zeta^{-j}d(d+1)\zeta^j=96(d - 3)(d - 1)^2(d - \frac12),      
    \end{equation*}
     is always non-zero, so there are two points of intersection. Moreover, the polynomials are independent of $k$, thus the intersection is independent of $k$, and the result follows. 
     
     The coordinates of the two points can be found by easy computation as \[\left(\frac{(d-2)(5d-3)\pm2(d-1)\sqrt{3(2d-1)(d-3)}}{d(d+1)}:\zeta^{-j}:0\right).\] 
   
   Note that when $d=3$, the discriminant above vanishes, and there is only one intersection point, $(1:\zeta^{-j}:0)$, and $V(z)$ is tangent to all three hyperosculating conics, as shown in \cite{DGMP25}.
\end{proof}

\begin{remark}\label{rem:twootherpoints}
    Note that between any two hyperosculating conics with fixed $j$, say $O_{j,k_1}$ and $O_{j,k_2}$, there are two other intersection points outside $V(z)$. Indeed, the intersection points between the conics are fixed in the pencil of the two, and the intersection points lie on the degenerate conic, $z\cdot \ell=0$, where (up to a constant) 
    \[\ell=2d(d-2)(\zeta^{-j}x+y)-2^{-1/d}(d+1)(2d-3)(u^{k_1}+u^{k_2})z.\] 
    The intersection of $O_{j,k_i}$ and $V(\ell)$ give the remaining two intersection points, and by inspection of the discriminant there are no values of $d \geq 3$ such that the two intersection points coalesce.
\end{remark}

\begin{remark}
 In \cite{DGMP25}, Dimca et. al. showed that the cubic Fermat curve together with three hyperosculating conics at three sextactic points co-linear with respect to $\mathcal{B}$, form a free curve. We suspect that this is not true for higher degrees.   
\end{remark}

We now consider sets of hyperosculating conics at sextactic points co-linear on the lines $\mathcal{M}$ and $\mathcal{N}$.

\begin{theorem}  
     When $d\geq 3$, the $d$ hyperosculating conics at $d$ sextactic points co-linear with respect to one of the lines in $\mathcal{M}$ (or $\mathcal{N}$) have two common intersection points.
\end{theorem}

\begin{proof}
    For a given $k$ odd, let $s_{j,k}$ for $j \in \{0,\ldots, d-1\}$ be the $d$ sextactic points on the line $z-u^{-k}2^{1/d}y=0$ contained in $\mathcal{M}_x$. 
     
     As before, we may use the defining polynomials from \Cref{cor:hyposccon} to compute the intersection points. Let $O_{j,k}$ be the $d$ hyperosculating conics at $s_{j,k}$. 
     
Observe that the $d$ points $s_{j,k}$ are in the same orbit under $\rho$, which fixes $V(x)$, so that the $d$ hyperosculating conics $O_{j,k}$ intersect $V(x)$ in the same points. 

The intersection of $V(x)$ with $O_{j,k}$ then gives  
     \begin{align*}
       V(d(d+1)u^{-k}2^{1/d}y^2+8d(d-2)yz-4(d+1)(2d-3)u^{k}2^{-1/d}z^2,x).
    \end{align*}
This intersection is independent of $j$, and the discriminant 
\begin{equation*}
\resizebox{\textwidth}{!}{$
    \left(8d(d-2)\right)^2-4\cdot d(d+1)u^{-k}2^{1/d}\cdot \left(-4(d+1)(2d-3)u^{k}2^{-1/d}\right)=48d(2d - 1)(d - 1)^2$}
\end{equation*}
is non-zero for all $d \geq 3$, so it is the same two intersection points for all the $d$ hyperosculating conics, 
\[\left(0:\frac{-4d(d-2)\pm2(d-1)\sqrt{3d(2d-1)}}{d(d+1)}:u^{-k}2^{1/d}\right).\] 
\end{proof}

\begin{remark}
    As in \Cref {rem:twootherpoints}, it is possible to show that there are two other intersection points (outside $V(x)$) between any two hyperosculating conics co-linear with respect to $\mathcal{M}$.
\end{remark}

\subsection{Intersections of osculating curves of degree \texorpdfstring{$n$}{n}}
Perhaps not surprisingly, the intersections of inflection tangents, tangents at sextactic points and hyperosculating conics to the Fermat curve can be generalized to intersections of osculating curves of any degree $n$. Since an automorphism preserves intersection multiplicities, we end this paper with an immediate corollary of \Cref{thm:invariantintersection}.
\begin{corollary}\label{cor:fullgenerality}
 Assume that $g$ is an automorphism of $\mathcal{F}$ with a line $L_{g}$ fixed pointwise by $g$. Let $P_{g}$ be the set of points in an orbit under $g$. Then the osculating curves of degree $n$ at the points in $P_g$ intersect $L_g$ in the same up to $n$ points, with the same respective intersection multiplicities at each intersection point.
\end{corollary}

\begin{remark}
   \Cref{cor:fullgenerality} can of course be generalized to any curve $\mathcal{C}$. The result is most interesting for the Fermat curves because of the abundance of lines fixed pointwise by automorphisms. 
\end{remark}

\begin{remark}
    Note that computing the osculating curve of degree $n>2$ to a point $p$ is not straightforward when $\mathcal{C}$ is not a rational curve.
\end{remark}

\bibliographystyle{plain}  
\bibliography{bib}    

\end{document}